\documentclass[english]{elsarticle}
\usepackage[T1]{fontenc}
\usepackage[latin9]{inputenc}
\usepackage{amsmath}
\usepackage{amssymb}
\usepackage{graphicx}
\usepackage{babel}
\begin{document}
\title{A low-rank algorithm for weakly compressible flow\tnoteref{label1}}

\author[tub,uibk]{\corref{cor1}Lukas Einkemmer}
\ead{lukas.einkemmer@uibk.ac.at}
\cortext[cor1]{Corresponding author}

\address[tub]{Department of Mathematics, University of T\"ubingen, Germany}
\address[uibk]{Department of Mathematics, University of Innsbruck, Austria}

\begin{abstract}
In this paper, we propose a numerical method for solving weakly compressible fluid flow based on a dynamical low-rank projector splitting. The low-rank splitting scheme is applied to the Boltzmann equation with BGK collision term, which results in a set of constant coefficient advection equations. This procedure is numerically efficient as a small rank is sufficient to obtain the relevant dynamics (described by the Navier--Stokes equations). The resulting method can be combined with a range of different discretization strategies; in particular, it is possible to implement spectral and semi-Lagrangian methods, which allows us to design numerical schemes that are not encumbered by the sonic CFL condition.

\end{abstract}
\maketitle

\section{Introduction}

Fluids play a pivotal role in virtually all fields of science and
engineering. Consequently, computational fluid dynamics (CFD) is used
from modeling pipe flows on a single workstation to simulating airplanes
or turbulent combustion on state of the art supercomputers. The governing
partial differential equations (PDEs) are the Navier\textendash Stokes
equations. More specifically, in the present work we will consider
the compressible isothermal Navier\textendash Stokes equations
\begin{align}
\partial_{t}\rho+\nabla_{x}(\rho u) & =0,\nonumber \\
\partial_{t}(\rho u)+\nabla\cdot(\rho u\otimes u)+\nabla p & =\nabla\cdot\left[\mu\left(\nabla u+(\nabla u)^{\mathrm{T}}\right)+\lambda(\nabla\cdot u)I\right]\label{eq:navier-stokes}\\
p & =\rho\theta,\nonumber 
\end{align}
where the density $\rho$ and the momentum $\rho u$ are the sought-after
quantities. Since we consider the isothermal case the (thermodynamic)
temperature $\theta$ is fixed. The pressure is determined by the
ideal gas law $p=\rho\theta$. Two material parameters, the dynamic
viscosity $\mu$ and the volume viscosity $\lambda$ have to be specified.
In the case of vanishing viscosity (i.e.~$\mu=0$ and $\lambda=0$)
equations (\ref{eq:navier-stokes}) are usually referred to as the
Euler equations.

The most common approach to solving these equations numerically is
to discretize them on an appropriate grid. Historically finite difference
and finite volume methods have been used extensively, while in recent
years discontinuous Galerkin schemes have become more common. However,
especially in the study of turbulence by direct numerical simulation
(DNS), spectral methods are often preferred (see, for example, \cite{kim1987turbulence,yokokawa200216}).

This approach (which we will refer to as direct discretization in
the following) is very mature and sophisticated numerical methods
have been developed in the last decades. Further advantages of this
approach are that (at least the basic) numerical algorithms are often
easy to understand and implement. Disadvantages include that explicit
methods usually need to satisfy the CFL condition for sound waves
(which in the weakly compressible setting can be multiple orders of
magnitude faster than the speed of flow) and that the equations are
relatively complicated (which puts significant constraints on the
design of numerical methods).

However, a direct discretization of the Navier\textendash Stokes equations
is not the only way to perform fluid simulations. In particular, lattice
Boltzmann, methods have been considered extensively in the literature
(see, for example, \cite{chen1998lattice,lallemand2000,kraemer2017}).
The lattice Boltzmann method exploits the fact the Boltzmann equation
(a kinetic model), for an appropriately modeled collision term and
initial value, recovers the dynamics of the Navier\textendash Stokes
equations (see, for example, \cite{bardos1991fluid,bardos1993fluid}).
Thus, in principle, we can solve fluid flow problems by integrating
the Boltzmann equation in time. However, these kinetic problems are
posed in a $2d$ dimensional phase space ($d$ dimensions of space,
as for the Navier\textendash Stokes equations, and $d$ dimensions
of velocity). Thus, a direct discretization is prohibitively expensive
from a computational point of view. The sought-after quantity, a distribution
function or particle-density, is usually denoted by $f(t,x,v)$.

However, in the fluid regime (i.e.~for thermalized gases or liquids)
we know that the distribution in velocity space stays close to a Maxwell\textendash Boltzmann
distribution. That is,
\[
f(t,x,v)\approx\frac{\rho(t,x)}{(2\pi\theta)^{d/2}}\exp\left(-\tfrac{1}{2}(v-u(t,x))^{2}\right).
\]
What we actually want to approximate are the moments of $f$ (which
correspond to the macroscopic quantities of density $\rho$ and momentum
$\rho u$). These are quantities of interest in fluid simulations
(as opposed to the distribution function $f$). The idea of the lattice
Boltzmann method is to discretize the velocity space with only a small
number of discrete velocities $e_{j}\in\mathbb{R}^{d}$. Then, the
moments can be computed using a Gaussian-type quadrature

\begin{align}
\rho(t,x) & =\int f(t,x,v)\,\mathrm{d}v\approx\sum_{j}W_{j}f_{j}(t,x),\label{eq:lb-moments-rho}\\
\rho u(t,x) & =\int vf(t,x,v)\,\mathrm{d}v\approx\sum_{j}W_{j}e_{j}f_{j}(t,x),\label{eq:lb-moments-rhou}
\end{align}
where $f_{j}(t,x)\approx f_{j}(t,x,e_{j})$ and $W_{j}$ are quadrature
weights. We then only have to solve an evolution equation for the
(relatively) small number of $f_{j}$s (which are $d$-dimensional
functions of $x$). For the classic lattice Boltzmann method in two-dimensions,
the $e_{j}$ are chosen as the corners of a square and the zero vector.
Thus, we have $9$ $f_{j}$ in two-dimensions (this is referred to
as D2Q9). In three dimensions a variety of schemes have been considered
(for example, D3Q19 and D3Q27 with $19$ and $27$ $f_{j}$s, respectively).
If the length of the square/cube is $2h/\tau$, where $h$ is the
grid spacing and $\tau$ is the time step size, a numerical method
(operated with unit CFL number) can be implemented without discretizing
any differential operators. This is a consequence of the fact that
the Boltzmann equation is much simpler compared to the Navier\textendash Stokes
equations. A further advantage of the lattice Boltzmann method is
that it can usually be parallelized very efficiently. Disadvantages
include that the amount of memory needed is increased (compared to
a direct discretization of the Navier\textendash Stokes equations)
and that the method is most effective if simulations are conducted
using a unit CFL number. Attempts to overcome the latter limitation
have resulted in the developed of so-called off-lattice Boltzmann
methods (see, for example, \cite{min2011,fakhari2015,kraemer2017}).
However, according to \cite{kraemer2017} these methods can be computationally
expensive due to the high number of partial differential equations
that have to be solved.

In the present paper we propose an alternative approach to both a
direct discretization of the Navier\textendash Stokes equations and
to the lattice Boltzmann method. Similar to the lattice Boltzmann
method our scheme is applied to the Boltzmann equation. However, to
reduce the dimensionality of the problem (from $2d$ to $d$) we perform
a low-rank approximation. We then obtain evolution equations that
describe the dynamics of the Boltzmann equation constraint to the
corresponding low-rank manifold. To accomplish this the dynamical
low-rank splitting algorithm introduced in \cite{Lubich2014} is used.
This allows us to represent the evolution in velocity space in more
detail. In fact, we obtain evolution equations for functions that
depend on $x$, but not on $v$ (as in the lattice Boltzmann method).
However, we also obtain similar evolution equations for functions
that depend only on $v$, but not on $x$. 

In addition, the evolution equations obtained are still significantly
simpler compared to the Navier\textendash Stokes equations (essentially
we obtain a constant-coefficient advection with an inhomogeneity).
Thus, a range of space discretization strategies can be employed relatively
easily. In particular, we can use (true) spectral methods (as opposed
to the pseudo-spectral approach which is common for the direct discretization
of the Navier\textendash Stokes equations). Furthermore, it is possible
within this approach to construct a numerical method that can overcome
the CFL condition imposed by the speed of sound. This is particular
relevant for weakly compressible simulations. Another interesting
property of the projector-splitting integrator is that it mimics the
properties of the (continuous) Boltzmann equation as we approach the
limit of vanishing viscosity (i.e.~as we consider the limit that
yields the Euler equations from the Navier\textendash Stokes equations).

Let us note that low-rank approximations have been extensively used
in quantum mechanics. See, in particular, \cite{meyer90tmc,meyer09mqd}
for the MCTDH approach to molecular quantum dynamics in the chemical
physics literature and \cite{Lubich2008,lubich15tii,Conte2010} for
a computational mathematics point of view. Some uses of dynamical
low-rank approximation in areas outside quantum mechanics are described
in \cite{Nonnenmacher2008,jahnke2008dynamical,Mena2017,musharbash2015error}.
In a general mathematical setting, dynamical low-rank approximation
has been studied in \cite{Koch2007,Koch2010,lubich13dab}. A major
algorithmic advance for the time integration was achieved with the
projector-splitting methods first proposed in \cite{Lubich2014} for
matrix differential equations and then developed further for various
tensor formats in \cite{lubich15tii,Lubich,Haegeman2016,Kieri2016,Lubich2017}.
Low-rank approximations for computational plasma physics (i.e.~the
collisionless but magnetized Boltzmann equation) have been considered
in \cite{Kormanna,einkemmer2018low}. Note, however, these schemes
try to capture kinetic effects that occur far away from thermodynamic
equilibrium. This means that the Navier\textendash Stokes equations
(or any other model that considers only the moments of $f$) are not
applicable in this setting. 

The outline of this paper is as follows. First, we summarize how the
Boltzmann equation gives rise to the Navier\textendash Stokes equations
(section \ref{sec:Boltztmann-to-fluid}). Then we introduce the proposed
numerical algorithm (section \ref{sec:numerical-method}). In section
\ref{sec:algorithm-in-the-limit} we investigate the behavior of the
low-rank projector-splitting as the viscosity vanishes. Numerical
results are presented in section \ref{sec:numerical-results}. Finally,
we conclude in section \ref{sec:conclusion}.

\section{Obtaining fluid dynamics from the Boltzmann equation\label{sec:Boltztmann-to-fluid}}

The compressible isothermal Navier\textendash Stokes equations have
already been stated in (\ref{eq:navier-stokes}). The goal in this
section is to show how the dynamics of the Navier\textendash Stokes
equations arises from the Boltzmann equation. Although, this has been
investigated before \cite{bardos1991fluid,bardos1993fluid}, it is
essential to motivate and explain the numerical algorithm that is
described in section \ref{sec:numerical-method}. 

In the remainder of the paper we will non-dimensionalize the Navier\textendash Stokes
equations as follows. First, we choose a characteristic length scale
$L$ (in section \ref{sec:numerical-results} this will be the length
of the computational domain). Then we choose the speed of sound $c_{s}$
as the characteristic velocity. Since the speed of sound is given
by $c_{s}=\sqrt{\frac{\partial p}{\partial\rho}}=\sqrt{\theta}$ this
choice implicitly sets $\theta=1$. As a consequence the characteristic
time is then $T=L/c_{s}$. This is just the time it takes a sound
wave to propagate from one end of the domain to the other end.

Since we focus on weakly compressible flow here, the flow speed is
always appreciably smaller than $1$ (the speed of sound) and we initialize
our problem with a homogeneous fluid density. The Navier\textendash Stokes
equations remain invariant with respect to scaling the density (i.e.~only
variations in the density are important). Thus, we simply initialize
$\rho=1$. For a real fluid this would determine the units used to
measure mass and thus also (uniquely) determines the units used to
measure viscosity. A flow with speed $u$ then has a Reynolds number
(in non-dimensionalized units)
\[
\text{Re}=\frac{\rho uL}{\mu}=\frac{u}{\mu}
\]
and a Mach number
\[
\text{M}=\frac{u}{c_{s}}=u.
\]

It is also instructive to discuss the incompressible limit. In this
case the (now incompressible) Navier\textendash Stokes equations become
\begin{align*}
\partial_{t}(\rho u)+\nabla\cdot(\rho u\otimes u)+\nabla p & =\mu\Delta u,\\
\nabla\cdot u & =0.
\end{align*}
Formally, this is obtained by taking $c_{s}\to\infty$. This also
explains the infinite speed of propagation modeled by the divergence
free constraint $\nabla\cdot u=0$. Clearly it is then not possible
to use the speed of sound as a characteristic velocity. Instead a
typical velocity of the flow is usually chosen. Care has to be taken
when comparing weakly compressible simulation results (such as those
in section \ref{sec:numerical-results}) with incompressible simulations
(such as those in \cite{bell1989second,liu2000high,einkemmer2014}).
In particular, in the latter the characteristic time is $T_{\text{incompr}}=L/u$
and in the former $T_{\text{compr}}=L/c_{s}=\text{M}L/u$. Thus, the
final time of a simulation has to be adjusted accordingly.

We now consider the Boltzmann equation
\begin{equation}
\partial_{t}f_{\epsilon}(t,x,v)+v\cdot\nabla_{x}f_{\epsilon}(t,x,v)=\frac{1}{\epsilon}C(f_{\epsilon})(x,v),\label{eq:boltzmann-generic-C}
\end{equation}
where $C$ is the collision operator and $\epsilon$ is a (usually
small) parameter. The sought-after quantity is the phase space distribution
$f_{\epsilon}.$ From a physical point of view the collision operator
has to enforce that the dynamics stays close to a Maxwell-{}-Boltzmann
distribution in velocity space. Consequently, we assume that
\begin{equation}
C(g)=0\label{eq:collision=00003D0}
\end{equation}
has solutions that can be written in the following form
\begin{equation}
g(x,v)=\frac{\rho(x)}{(2\pi\theta)^{d/2}}\exp\left(-\frac{1}{2}\frac{(v-u(x))^{2}}{\theta}\right),\label{eq:equi-C}
\end{equation}
where the density $\rho$ and the velocity $u$ (strictly speaking,
the momentum $\rho u$), are given by the moments
\[
\rho(x)=\int g(x,v)\,\mathrm{d}v,\qquad\qquad\rho u(x)=\int vg(x,v)\,\mathrm{d}v.
\]
From now on we set the (thermodynamic) temperature to one (i.e.~$\theta=1$).
Equation (\ref{eq:equi-C}) is precisely what we would, on physical
grounds, expect from an ideal thermalized fluid. 

The remarkable observation here is that, equation (\ref{eq:boltzmann-generic-C})
still fully captures the (very complicated) dynamics of the Navier\textendash Stokes
equations. Conceptually the simplest case is the limit $\epsilon\to0$.
Thus, we will consider it here. In this case the right-hand side of
equation (\ref{eq:boltzmann-generic-C}) constrains the solution to
the form
\begin{equation}
f^{\text{eq}}(t,x,v)=\frac{\rho(t,x)}{(2\pi)^{d/2}}\exp\left(-\tfrac{1}{2}(v-u(t,x))^{2}\right),\label{eq:f-dyn-equi}
\end{equation}
where $\rho$ and $u$ are, as of yet, undetermined quantities. Thus,
we have $f^{\text{eq}}=\lim_{\epsilon\to0}f_{\epsilon}$. We proceed
by integrating equation (\ref{eq:boltzmann-generic-C}) with respect
to velocity and obtain
\[
\partial_{t}\int f_{\epsilon}\,\mathrm{d}v+\nabla_{x}\cdot\left(\int vf_{\epsilon}\,\mathrm{d}v\right)=\frac{1}{\epsilon}C(f_{\epsilon})
\]
Now, we take the (formal) limit $\epsilon\to0$ 

\[
\partial_{t}\int f^{\text{eq}}(t,x,v)\,\mathrm{d}v+\nabla_{x}\cdot\left(\int vf^{\text{eq}}(t,x,v)\,\mathrm{d}v\right)=0.
\]
Note that the collision term has vanished as the solution $f^{\text{eq}}$
given by equation (\ref{eq:f-dyn-equi}) satisfies (\ref{eq:collision=00003D0}).
Employing the definition of density and velocity as the moments of
the phase space density, we easily obtain
\[
\partial_{t}\rho+\nabla_{x}\cdot(\rho u)=0.
\]
This is precisely the continuity equation.

To derive the momentum balance equations, we multiply equation (\ref{eq:boltzmann-generic-C})
by $v_{j}$ (the $j$th component of the velocity) and integrate in
velocity space. This yields
\[
\partial_{t}\int v_{j}f_{\epsilon}\,\mathrm{d}v+\sum_{i}\partial_{x_{i}}\int v_{j}v_{i}f_{\epsilon}\,\mathrm{d}v=\frac{v_{j}}{\epsilon}C(f_{\epsilon}).
\]
Taking the limit $\epsilon\to0$ gives
\[
\partial_{t}(\rho u_{j})+\sum_{i}\partial_{x_{i}}\int v_{j}v_{i}f^{\text{eq}}\,\mathrm{d}v=0.
\]
We now evaluate the resulting integrals by using equation (\ref{eq:f-dyn-equi}).
We obtain
\[
\int v_{j}v_{i}f^{\text{eq}}\,\mathrm{d}v=\begin{cases}
\rho\left(u_{i}^{2}+1\right) & i=j\\
\rho u_{i}u_{j} & i\neq j
\end{cases}
\]
and thus
\[
\partial_{t}(\rho u)+\nabla\cdot(\rho u\otimes u)+\nabla\rho=0.
\]
This is precisely the momentum balance equation for an ideal gas.
More canonically we would write this using the pressure $p$ and impose
the equation of state corresponding to an ideal gas, i.e.~$p=\rho$.
In summary, we have obtained the compressible isothermal Euler equations.

The question that remains to be answered is why a low-rank representation
makes sense here. We know that the solution satisfies the form specified
by equation (\ref{eq:f-dyn-equi}) at all times. However, this is
not a low-rank representation due to the presence of both velocity
($v$) and position ($x$) dependent functions in the exponential.
However, if the flow velocity is small compared to the speed of sound
(i.e.~in the weakly compressible case) we can use
\begin{equation}
f^{\text{eq}}=\frac{\rho}{(2\pi)^{d/2}}\exp\left(-\frac{v^{2}}{2}\right)\left(1+v\cdot u+\frac{(v\cdot u)^{2}}{2}-\frac{u^{2}}{2}\right)+\mathcal{O}\left(u^{3}\right).\label{eq:feq-expansion}
\end{equation}
This is a low-rank approximation with rank $6$ and $10$ for two-
and three-dimensional problems, respectively. For comparison, a lattice
Boltzmann method usually requires $9$ directions in two-dimensions
and $19$ to $27$ directions in three-dimensions (see the discussion
in the introduction). Thus, at least in princple, representing the
solution by a low-rank repersentation is a viable approach.

The derivation for the Navier\textendash Stokes equations (i.e.~for
$\epsilon>0$) is more involved. One proceeds by performing a Chapman\textendash Enskog
expansion. That is, we assume that $\epsilon$ is a small parameter
and look for a solution, up to terms of $\mathcal{O}\left(\epsilon^{2}\right)$,
to equation (\ref{eq:boltzmann-generic-C}) that has the form
\begin{equation}
f_{\epsilon}=f^{\text{eq}}(1+\epsilon g_{\epsilon}+\epsilon^{2}w_{\epsilon}),\label{eq:feps-expansion-1}
\end{equation}
where $f^{\text{eq}}$ is given, as before, in the form specified
by equation (\ref{eq:f-dyn-equi}). The functions $g_{\epsilon}$
and $w_{\epsilon}$ give, respectively, the first and second order
deviation from $f^{\text{eq}}$ caused by the finite $\epsilon$.
We will use the BGK (Bhatnagar\textendash Gross\textendash Krook)
collision operator
\[
C(f_{\epsilon})=f^{\text{eq }}-f_{\epsilon}.
\]
It is easy to check that this collision operator satisfies the condition
given in (\ref{eq:collision=00003D0}). The BGK collision operator
is heavily used in lattice Boltzmann simulations and we will also
employ it for the numerical results conducted in section \ref{sec:numerical-results}.

It can then be shown that we recover the continuity equation
\[
\partial_{t}\rho+\nabla\cdot(\rho u)=0
\]
and the following momentum balance equations
\begin{equation}
\partial_{t}(\rho u)+\nabla\cdot(\rho u\otimes u)+\nabla\rho=\epsilon\nabla\cdot\left[\rho\nabla u+\rho(\nabla u)^{\text{T}}-\tfrac{2}{d}\rho(\nabla\cdot u)I\right].\label{eq:derived-ns-raw-1}
\end{equation}
For more details of the derivation we refer the reader to \cite{bardos1991fluid,bardos1993fluid,xu2003lattice}.
Comparing this to the Navier\textendash Stokes equations stated in
the introduction, i.e.~equation (\ref{eq:navier-stokes}), we have
a perfect match, except for the diffusion term. However, in the case
of weakly compressible flow $\rho$ varies only slightly. In addition,
non-dimensionalization allows us to set the characteristic value of
$\rho$ to unity. Thus, making the approximation $\rho\approx1$ and
applying it to the right-hand side of equation (\ref{eq:derived-ns-raw-1})
we obtain
\[
\partial_{t}(\rho u)+\nabla\cdot(\rho u\otimes u)+\nabla\rho=\epsilon\nabla\cdot\left[\nabla u+(\nabla u)^{\text{T}}-\frac{2}{d}(\nabla\cdot u)I\right].
\]
Thus, we have recovered the Navier\textendash Stokes equation with
$\mu=\epsilon$ and $\lambda=-2\epsilon/d$.

It should also be noted that, while there is no guarantee, that $g_{\epsilon}$
is a low-rank function, the dynamics implied by the Navier\textendash Stokes
equations forces $f_{\epsilon}$ to stay close to a low-rank function.
This further motivates the proposed approach and we will see in section
\ref{sec:numerical-results} that usually quite low ranks are sufficient
in order to obtain excellent agreement with the dynamics of interest.

\section{Numerical method\label{sec:numerical-method}}

We start from the Boltzmann equation
\begin{equation}
\partial_{t}f(t,x,v)+v\cdot\nabla_{x}f(t,x,v)=\frac{1}{\epsilon}C(f)(x,v),\label{eq:boltzmann-bkg}
\end{equation}
with the BGK collision operator
\[
C(f)=f^{\text{eq}}-f,
\]
where $\epsilon>0$ is a (usually small) parameter and
\[
f^{\text{eq}}(t,x,v)=\frac{\rho(t,x)}{(2\pi)^{d/2}}\exp\left(-\tfrac{1}{2}(v-u(t,x))^{2}\right),
\]
where $d\in\{1,2,3\}$ is the dimension of the problem. The sought-after
quantity is $f$ (in this section we will not explicitly denote the
dependence of $f$ on $\epsilon$). As has been outlined in the previous
section, the moments 
\[
\rho=\int f\,\mathrm{d}v,\qquad\qquad\rho u=\int vf\,\mathrm{d}v
\]
then satisfy the compressible isothermal Navier\textendash Stokes
equations. As initial value we choose a function of the form
\[
f(0,x,v)=\frac{\rho^{0}(x)}{(2\pi)^{d/2}}\exp\left(-\tfrac{1}{2}\left(v-u^{0}(x)\right)^{2}\right),
\]
This is not yet a low-rank representation. However, in an actual implementation
we can either use the expansion given in equation (\ref{eq:feq-expansion})
or perform a singular value decomposition (SVD) once the problem is
discretized.

What remains to be determined here is the density $\rho^{0}$ and
the velocity $u^{0}$ (or alternatively, the momentum $\rho u^{0}$).
These are specified according to the fluid problem for which a numerical
solution is sought.

Since equation (\ref{eq:boltzmann-bkg}) is posed in a $2d$ dimensional
phase space, its direct solution is prohibitively expensive. This
is particularly true in the present setting as the dynamics stays
close to a low-rank manifold (see the discussion in the previous section).
Thus, the goal of this section is to derive an algorithm that approximates
the Boltzmann equation (\ref{eq:boltzmann-bkg}) by a low-rank representation.

To that end, the function $f(t,x,v)$ is constrained to the following
form
\begin{equation}
f(t,x,v)=\sum_{ij}X_{i}(t,x)S_{ij}(t)V_{j}(t,v),\label{eq:f-lowrank}
\end{equation}
where $S\in\mathbb{R}^{r\times r}$ and we call $r$ the rank of the
representation. Note that the dependence of $f$ on the phase space
$(x,v)\in\Omega\subset\mathbb{R}^{2d}$ is now approximated by the
functions $\{X_{i}\colon i=1,\dots,r\}$ and $\{V_{j}\colon j=1,\dots,r\}$
which depend only on $x\in\Omega_{x}\subset\mathbb{R}^{d}$ and $v\in\Omega_{v}\subset\mathbb{R}^{d}$,
respectively. In equation (\ref{eq:f-lowrank}) and the following
discussion we always assume that summation indices run from $1$ to
$r$ and we thus do not specify these bounds.

Now, we seek an approximation to the exact particle-density function
that for all $t$ lies in the set
\[
\overline{\mathcal{M}}=\biggl\{ f\in L^{2}(\Omega)\colon f(x,v)=\sum_{ij}X_{i}(x)S_{ij}V_{j}(v)\text{ with }S\in\mathbb{R}^{r\times r},\,X_{i}\in L^{2}(\Omega_{x}),\,V_{j}\in L^{2}(\Omega_{v})\biggr\}.
\]
It is clear that this representation is not unique. In particular,
we can make the assumption that $(X_{i},X_{k})=\delta_{ik}$ and $(V_{j},V_{l})=\delta_{jl}$,
where $(\cdot,\cdot)$ is the inner product on $L^{2}(\Omega_{x})$
and $L^{2}(\Omega_{v})$, respectively. We consider a path $f(t)$
on $\overline{\mathcal{M}}$. The corresponding derivative is denoted
by $\dot{f}$ and is of the form 
\begin{equation}
\dot{f}=\sum_{ij}\left(X_{i}\dot{S}_{ij}V_{j}+\dot{X}_{i}S_{ij}V_{j}+X_{i}S_{ij}\dot{V}_{j}\right).\label{eq:path}
\end{equation}
If we impose the conditions $(X_{i},\dot{X}_{j})=(V_{i},\dot{V}_{j})=0$
then $S_{ij}$ is uniquely determined by $\dot{f}$. This follows
easily from the fact that
\begin{equation}
\dot{S}_{ij}=(X_{i}V_{j},\dot{f}).\label{eq:unique-S}
\end{equation}
We then project both sides of equation (\ref{eq:path}) onto $X_{i}$
and $V_{j}$, respectively, and obtain
\begin{align}
\sum_{j}S_{ij}\dot{V}_{j} & =(X_{i},\dot{f})-\sum_{j}\dot{S}_{ij}V_{j},\label{eq:unique-V}\\
\sum_{i}S_{ij}\dot{X}_{i} & =(V_{j},\dot{f})-\sum_{i}X_{i}\dot{S}_{ij}.\label{eq:unique-X}
\end{align}
From these relation it follows that the $X_{i}$ and $V_{j}$ are
uniquely defined if $S$ has full rank (this, in particular, implies
that $S$ and $S^{T}$ are invertible). Thus, we seek an approximation
that for each time $t$ lies in the manifold
\begin{align*}
\mathcal{M} & =\biggl\{ f\in L^{2}(\Omega)\colon f(x,v)=\sum_{ij}X_{i}(x)S_{ij}V_{j}(v)\text{ with }S\in\mathbb{R}^{r\times r},\,X_{i}\in L^{2}(\Omega_{x}),\,V_{j}\in L^{2}(\Omega_{v})\text{ and}\\
 & \qquad\qquad(X_{i},X_{k})=\delta_{ik},\,(V_{j},V_{l})=\delta_{jl},S\text{ has full rank}\biggr\}
\end{align*}
with the corresponding tangent space 
\begin{align*}
\mathcal{T}_{f}\mathcal{M} & =\biggl\{\dot{f}\in L^{2}(\Omega)\colon\dot{f}(x,v)=\sum_{ij}\left(X_{i}(x)\dot{S}_{ij}V_{j}(v)+\dot{X}_{i}(x)S_{ij}V_{j}(v)+X_{i}(x)S_{ij}\dot{V}_{j}(v)\right),\\
 & \qquad\qquad\text{with }\dot{S}\in\mathbb{R}^{r\times r},\,\dot{X}_{i}\in L^{2}(\Omega_{x}),\,\dot{V}_{j}\in L^{2}(\Omega_{v}),\text{ and }(X_{i},\dot{X}_{j})=(V_{i},\dot{V_{j}})=0\biggr\},
\end{align*}
where $f$ is given by equation (\ref{eq:f-lowrank}). Now, we consider
the dynamics of the Boltzmann equation on the manifold $\mathcal{M}$.
That is, we consider
\begin{equation}
\partial_{t}f=-P(f)\left(v\cdot\nabla_{x}f-\frac{1}{\epsilon}C(f)\right),\label{eq:vlasov-proj}
\end{equation}
where $P(f)$ is the orthogonal projector onto the tangent space $\mathcal{T}_{f}\mathcal{M}$,
as defined above.

We will consider the projection $P(f)g$ for a moment. From equations
(\ref{eq:path})-(\ref{eq:unique-X}) we obtain
\[
P(f)g=\sum_{j}(V_{j},g)V_{j}-\sum_{ij}X_{i}(X_{i}V_{j},g)V_{j}+\sum_{i}(X_{i},g)X_{i}.
\]
Let us introduce the following two vector spaces $\overline{X}=\text{span}\left\{ X_{i}\colon i=1,\dots,r\right\} $
and $\overline{V}=\text{span}\left\{ V_{j}\colon j=1,\dots r\right\} $.
Then we can write the projector as follows

\begin{equation}
P(f)g=P_{\overline{V}}g-P_{\overline{V}}P_{\overline{X}}g+P_{\overline{X}}g,\label{eq:projector}
\end{equation}
where $P_{L}$ is the orthogonal projector onto the vector space $L$.
The decomposition of the projector into this three terms forms the
basis of our splitting procedure (for matrix equations this has been
first suggested in \cite{Lubich2014}). 

We proceed by substituting $g=v\cdot\nabla_{x}f-\frac{1}{\epsilon}C(f)$
into equation (\ref{eq:projector}). This at once gives a three-term
splitting for equation (\ref{eq:vlasov-proj}). More precisely, for
the first order Lie splitting we have to solve the equations
\begin{align}
\partial_{t}f & =-P_{\overline{V}}\left(v\cdot\nabla_{x}f-\frac{1}{\epsilon}C(f)\right),\label{eq:split-I}\\
\partial_{t}f & =P_{\overline{V}}P_{\overline{X}}\left(v\cdot\nabla_{x}f-\frac{1}{\epsilon}C(f)\right)\label{eq:split-II}\\
\partial_{t}f & =-P_{\overline{X}}\left(v\cdot\nabla_{x}f-\frac{1}{\epsilon}C(f)\right)\label{eq:split-III}
\end{align}
one after another. In the following discussion we will consider the
first order Lie splitting algorithm with step size $\tau$.

We assume that the initial value for the algorithm is given in the
following form
\[
f(0,x,v)=\sum_{ij}X_{i}^{0}(x)S_{ij}^{0}V_{j}^{0}(v).
\]

First, let us consider equation (\ref{eq:split-I}). Since the set
$\{V_{j}\colon j=1,\dots,r\}$ forms an orthonormal basis of $\overline{V}$
(for each $t$), we have
\begin{equation}
f(t,x,v)=\sum_{j}K_{j}(t,x)V_{j}(t,v),\qquad\quad K_{j}(t,x)=\sum_{i}X_{i}(t,x)S_{ij}(t),\label{eq:Kdef}
\end{equation}
where $K_{j}(t,x)$ is the coefficient of $V_{j}$ in the corresponding
basis expansion. We duly note that $K_{j}$ is a function of $x$,
but not of $v$). We then rewrite equation (\ref{eq:split-I}) as
follows
\begin{align*}
 & \sum_{j}\partial_{t}K_{j}(t,x)V_{j}(t,v)+\sum_{j}K_{j}(t,x)\partial_{t}V_{j}(t,v)\\
 & =-\sum_{j}\left(V_{j}(t,\cdot),v\mapsto v\cdot\nabla_{x}f(t,x,v)-\tfrac{1}{\epsilon}C(f)(x,v))\right)V_{j}(t,v).
\end{align*}
The solution of this equation is given by $V_{j}(t,v)=V_{j}(0,v)=V^{0}(v)$
and
\begin{equation}
-\dot{K}_{j}(t,x)=\sum_{l}c_{jl}^{1}\cdot\nabla_{x}K_{l}(t,x)-\frac{1}{\epsilon}\left(K_{j}-c_{j}^{3}(K)(x)\rho(K)(x)\right)\label{eq:evolution-K}
\end{equation}
with
\begin{align*}
c_{jl}^{1} & =\int vV_{j}^{0}V_{l}^{0}\,\mathrm{d}v,\qquad c_{j}^{3}(K)(x)=\int V_{j}h^{\text{eq}}(K)\,\mathrm{d}v,
\end{align*}
where we have used the decomposition $f^{\text{eq}}=\rho h^{\text{eq}}$.
The evolution equation is obtained by equating coefficients in the
basis expansion. A very useful property of the present splitting is
that we have to only update the $K_{j}$, but not the $V_{j}$. We
further note that $c_{jl}^{1}=(c_{jl}^{1;x_{1}},c_{jl}^{1;x_{2}})$
(for $d=2$) is a vector quantity. Also note that we use $c_{j}^{3}$
here (instead of $c_{j}^{2}$) to keep the notation in line with \cite{einkemmer2018low},
where $c_{j}^{2}$ was used for the term originating from the electric
field (which is not present for standard fluid flow). However, since,
as is briefly discussed in the conclusion, the proposed numerical
method could conceivably be generalized to magnetohydrodynamic problems,
we have chosen this notation. 

Equation (\ref{eq:evolution-K}) is completely posed in a $d$-dimensional
(as opposed to $2d$-dimensional) space. Thus, we proceed by integrating
equation (\ref{eq:evolution-K}) with initial value
\[
K_{j}(0,x)=\sum_{i}X_{i}^{0}(x)S_{ij}^{0}
\]
until time $\tau$ to obtain $K_{j}^{1}(x)=K_{j}(\tau,x)$. However,
this is not sufficient as the $K_{j}^{1}$ are not necessarily orthogonal
(a requirement of our low-rank representation). Fortunately, this
is easily remedied by performing a QR decomposition
\[
K_{j}^{1}(x)=\sum_{i}X_{i}^{1}(x)S_{ij}^{1}
\]
to obtain orthonormal $X_{i}^{1}$ and the matrix $S_{ij}^{1}$. Once
a space discretization has been introduced, this QR decomposition
can be simply computed by using an appropriate function from a software
package such as LAPACK. However, from a mathematical point of view,
the continuous dependence on $x$ causes no issues. For example, the
modified Gram-Schmidt process works just as well in the continuous
formulation considered here.

Second, we proceed in a similar way for equation (\ref{eq:split-II}).
In this case both $V_{j}^{0}$ and $X_{i}^{1}$ are unchanged and
only $S_{ij}$ is updated. The corresponding evolution equation is
given by
\begin{align}
\dot{S}_{ij}(t) & =\left(X_{i}^{1}V_{j}^{0},(x,v)\mapsto(v\cdot\nabla_{x}f(t,x,v)-\tfrac{1}{\epsilon}C(f)(x,v))\right)\nonumber \\
 & =\sum_{lk}(d_{il}^{1}\cdot c_{jk}^{1})S_{lk}+\frac{1}{\epsilon}\left(S_{ij}(t)-e_{ij}(S)\right)\label{eq:evolution-S}
\end{align}
with
\[
d_{il}^{1}=\int X_{i}^{1}\nabla_{x}X_{l}^{1}\,\mathrm{d}x,\qquad\quad e_{ij}=\int X_{i}^{1}\rho(S)V_{j}^{0}h^{\text{eq}}(S)\,\mathrm{d}(x,v).
\]
Note that in this case the evolution equation depends neither on $x$
nor on $v$. We now integrate equation (\ref{eq:evolution-S}) with
initial value $S_{ij}(0)=S_{ij}^{1}$ until time $\tau$ and obtain
$S_{ij}^{2}=S_{ij}(\tau)$. This completes the second step of the
algorithm.

Finally, we consider equation (\ref{eq:split-III}). Similar to the
first step we have
\[
f(t,x,v)=\sum_{i}X_{i}(t,x)L_{i}(t,v),\qquad\quad L_{i}(t,v)=\sum_{j}S_{ij}(t)V_{j}(t,v).
\]
As before, it is easy to show that the $X_{i}$ remain constant during
that step. Thus, the $L_{j}$ satisfy the following evolution equation
\begin{align}
-\dot{L}_{i}(t,v) & =\left(X_{j}^{1},x\mapsto(v\cdot\nabla_{x}f(t,x,v)-\tfrac{1}{\epsilon}C(f)(x,v))\right)\nonumber \\
 & =-\sum_{l}(d_{il}^{1}\cdot v)L_{l}-\tfrac{1}{\epsilon}\left(L_{i}-d_{i}^{3}(L)(v)\right)\label{eq:evolution-L}
\end{align}
with
\[
d_{i}^{3}(L)(v)=\int X_{i}\rho h^{\text{eq}}\,dx.
\]
We then integrate equation (\ref{eq:evolution-L}) with initial value
\[
L_{i}(0,v)=\sum_{j}S_{ij}^{2}V_{j}^{0}(v)
\]
up to time $\tau$ to obtain $L_{i}^{1}(v)=L_{i}(\tau,v)$. Since,
in general, the $L_{i}^{1}$ are not orthogonal we have to perform
a QR decomposition
\[
L_{i}^{1}(v)=\sum_{j}S_{ij}^{3}V_{j}^{1}(v)
\]
to obtain $S_{ij}^{3}$ and $V_{j}^{1}$. Finally, the output of our
Lie splitting algorithm is
\[
f(\tau,x,v)\approx\sum_{ij}X_{i}^{1}(x)S_{ij}^{3}V_{j}^{1}(v).
\]
For simplicity, we have introduced the low-rank algorithm in the context
of the first order Lie splitting here. However, the extension to second-order
Strang splitting, which we use in the numerical simulations conducted
in section \ref{sec:numerical-results}, is straightforward.

Note that, to some extend, the algorithm introduced here has certain
similarities with a lattice Boltzmann method. In particular, the $X_{i}(t,x)$
roughly correspond to the $f_{i}(t,x)$ in the introduction. However,
there are important differences. In a lattice Boltzmann method the
distribution function $f$ would be represented as
\[
f(t,x,v)=\sum_{i}W_{i}f_{i}(t,x)\delta(v-e_{i}).
\]
This yields the correct moments according to equations (\ref{eq:lb-moments-rho})
and (\ref{eq:lb-moments-rhou}). For the proposed algorithm, however,
we consider the functions $V_{j}(t,v)$ which are propagated in time.
Thus, we do not consider only a single velocity per $X_{i}$ but rather
a distribution of velocities.

\subsection{Discretization}

The evolution equations for $S_{ij}$ and $L_{i}$ do not involve
any spatial derivatives and thus require no further discretization
(with the exception of the coefficients, which are constant during
the corresponding sub-step).

However, the evolution equation that describe the dynamics of the
$K_{j}$ are given by (for simplicity we only consider the two-dimensional
case here; however, the extension to three dimensions is immediate)
\begin{equation}
\partial_{t}K_{j}=-\sum_{l}(c_{jl}^{1;x_{1}}\partial_{x_{1}}K_{l}+c_{jl}^{1;x_{2}}\partial_{x_{2}}K_{l})-\frac{1}{\epsilon}\left(K_{j}-c_{j}^{3}(K)(x)\rho(K)(x)\right).\label{eq:evol-K-detail2d}
\end{equation}
Since this is a constant-coefficient advection, we can choose virtually
any space discretization scheme (finite differences, finite volumes,
etc.) to obtain 
\begin{equation}
\partial_{t}K=AK-\frac{1}{\epsilon}(K-c_{j}^{3}(K)\rho(K)),\label{eq:evol-K-matrix}
\end{equation}
where $K=[K_{1},\dots,K_{r}]$ and $A$ is a matrix that represents
the discretized differential operator.

Let us pause here for a moment. In the literature a number of different
techniques have been developed to solve the Euler equations (or, more
generally, fluid flow where sharp gradients occur). Often such techniques
are based on upwind schemes. While implementing upwind schemes for
a scalar constant coefficient advection equation is a rather simple
task, the (non-scalar and nonlinear) nature of the Euler equations
makes this significantly more challenging in practice. For a good
review we refer the reader to \cite{van2006upwind}. One way to generalize
upwind schemes is to solve a Riemann problem at the cell interface,
which can incur a significant computational cost. Now, note that since
equations (\ref{eq:evol-K-detail2d}) are constant-coefficient advections,
most of these difficulties are avoided for the numerical scheme proposed
here. Thus, upwind schemes can be implemented relatively easily as
part of the proposed numerical algorithm. We will not explore this
topic further in the present paper, but we consider this as future
work.

In principle, equation (\ref{eq:evol-K-matrix}) can be solved by
an appropriate time integrator. Note, however, that using an explicit
method would introduce a CFL condition. In this case we would use
sub-stepping. That is, a smaller time step is used to solve (\ref{eq:evol-K-matrix})
compared to the splitting scheme. However, this can be avoided by
employing a semi-Lagrangian approach (as discussed in the following)
or a spectral approach (as discussed in the next section). To do that
we first apply a further splitting procedure to equation (\ref{eq:evol-K-detail2d}).
For Lie splitting this yields
\[
K(\tau,\cdot)\approx\varphi_{\tau}^{\epsilon}\left(\mathrm{e}^{-\tau c^{1;x_{2}}\partial_{x_{2}}}\mathrm{e}^{-\tau c^{1;x_{1}}\partial_{x_{1}}}K(0,\cdot)\right),
\]
where $\varphi_{\tau}^{\epsilon}$ is the partial flow generated by
the collision term. The crucial part is the computation of
\[
M(t,x)=\mathrm{e}^{-\tau c^{1;x_{1}}\partial_{x_{1}}}M(0,x)
\]
which is equivalent to the partial differential equation
\[
\partial_{t}M(t,x)=-c^{1;x_{1}}\partial_{x_{1}}M(t,x).
\]
Now, since $c^{1;x_{1}}$ is symmetric, there exists an orthogonal
matrix $T$ such that $Tc^{1;x_{1}}T^{T}=D$, where $D$ is a diagonal
matrix. All the ingredients can be computed efficiently as $c^{1;x_{1}}\in\mathbb{R}^{r\times r}$
(i.e.~these are small matrices). We now change variables to $\overline{M}=TM$
and obtain
\[
\partial_{t}\overline{M}_{j}(t,x)=-D_{jj}\partial_{x_{1}}\overline{M}(t,x).
\]
This is now a set of scalar one-dimensional advection equation with
constant coefficients and can thus be treated by an arbitrary semi-Lagrangian
approach. 

\subsection{Spectral discretization}

Pseudo-spectral methods are widely used in some fluid problems (for
example, for turbulent DNS simulations \cite{kim1987turbulence,yokokawa200216}).
Here we will show that (true) spectral methods can be very naturally
incorporated into the proposed low-rank scheme. To do that we perform
the Fourier transformation with respect to $x$ of equation (\ref{eq:evol-K-detail2d}).
This yields
\begin{equation}
\partial_{t}\hat{K}(t,k)=A(k)\hat{K}(t,k)-\frac{1}{\epsilon}(\hat{K}-c^{3}(K)\rho(K)),\label{eq:spectral-discr}
\end{equation}
where $\hat{K}_{j}$ denotes the Fourier transform of $K_{j}$ and
we have defined $\hat{K}=[\hat{K}_{1},\dots,\hat{K}_{r}]$ and $K=[K_{1},\dots,K_{r}]$.
This would be sufficient for a pseudo-spectral approach. However,
we can turn this into a spectral method by further splitting equation
(\ref{eq:spectral-discr}). This is possible since we only have to
treat constant-coefficient advection equations and the nonlinear term
(i.e.~the collision operator) is free of spatial derivatives. In
particular, this is in contrast to the Navier\textendash Stokes equations,
where the nonlinear terms involve spatial differentiation. For Lie
splitting this yields
\[
\hat{K}(\tau,k)\approx\varphi_{\tau}^{\epsilon}\left(\mathrm{e}^{\tau A(k)}\hat{K}(0,k)\right),
\]
where $\varphi_{\tau}^{\epsilon}$ is the partial flow generated by
the collision term. The exponential can be readily computed in Fourier
space as $A(k)\in\mathbb{R}^{r\times r}$ (and thus we only have to
compute the exponential of a small matrix). We also note that this
approach is, obviously, not encumbered by a CFL condition.

\subsection{Computational efficiency}

In this section, we will discuss the computational characteristics
of the proposed algorithm. Solving the evolution equations is at most
$\mathcal{O}\left(rn^{d}\right)$ (both in terms of cost as well as
in terms of storage), where $n$ is the number of grid points per
direction. As we will see in section \ref{sec:numerical-results},
it is often sufficient to use significantly fewer grid points in the
velocity (i.e.~$v$) directions than in the spatial (i.e.~$x$)
directions. Thus, the evolution equation for $K$, equation (\ref{eq:evolution-K}),
dominates the computational effort. This, in particular, makes the
comparison to lattice Boltzmann methods (which only have to integrate
$x$-dependent quantities) more favorable.

However, in addition, for the proposed numerical method we have to
compute various coefficients. To compute the coefficients $c_{jl}^{1}$
and $d_{il}^{1}$ requires a computational cost of $\mathcal{O}\left(r^{2}n^{d}\right)$
and $\mathcal{O}\left(r^{2}\right)$ storage. Now, naively computing
$c_{j}^{3}$ and $d_{i}^{3}$ would be quite expensive and could easily
dominate the run time of our algorithm. However, we can accomplish
this with a computational cost of $\mathcal{O}\left(r^{2}n^{d}\right)$.
To do that we proceed as follow. First, we write 
\[
h^{\text{eq}}(x,v)=\frac{1}{(2\pi)^{d/2}}\exp\left(-\frac{v^{2}}{2}\right)\sum_{k}h_{k}^{X}(x)h_{k}^{V}(v),
\]
where the sum is over $10$/$6$ ($d=3$/$2$) entries and each $h^{X}$
and $h^{V}$ is a monomial (see the expansion in equation (\ref{eq:feq-expansion})).
Thus, we exploit the low-rank expansion of $h^{\text{eq}}$. Then
we rewrite $c_{j}^{3}$ as follows
\[
c_{j}^{3}(x)=\sum_{k}h_{k}^{X}(x)I_{jk}^{1},\qquad\qquad I_{jk}^{1}=\frac{1}{(2\pi)^{d/2}}\int V_{j}(v)\exp\left(-\frac{v^{2}}{2}\right)h_{k}^{V}(v)\,\mathrm{d}v
\]
and $d_{j}^{3}$ as follows
\[
d_{i}^{3}(v)=\frac{1}{(2\pi)^{d/2}}\exp\left(-\frac{v^{2}}{2}\right)\sum_{k}h_{k}^{V}(v)I_{ik}^{2},\qquad\qquad I_{ik}^{2}=\int X_{i}\rho h_{k}^{X}\,\mathrm{d}x.
\]
Both computing the integrals and summing the results to obtain $c_{j}^{3}$
and $d_{i}^{3}$ requires a computational cost of $\mathcal{O}\left(r^{2}n^{d}\right)$.
Finally, we can use $c_{j}^{3}$ to compute $e_{ij}$ as follows
\[
e_{ij}=\int X_{i}\rho c_{j}^{3}\,dx.
\]
This has a computational cost of $\mathcal{O}\left(r^{2}n^{d}\right)$.
Thus, the entire algorithm can be implemented with a computational
cost of $\mathcal{O}\left(r^{2}n^{d}\right)$ and a storage cost of
$\mathcal{O}\left(rn^{d}\right)$. In practice, computing these coefficients
might even be faster as, for example, the computation of $I_{ik}^{2}$
is limited by the memory loads of $X_{i}$ and $\rho$ ($h_{k}^{X}$
is a monomial which we can easily computed on the fly).

One might be worried that the proposed algorithm requires $\mathcal{O}\left(r^{2}n^{d}\right)$
arithmetic operations. However, we only require $\mathcal{O}\left(rn^{d}\right)$
memory operations. The latter, in a reasonable implementation, dominates
the performance of the algorithm on all present and, most likely,
all future computer systems. A hope is that (especially in three-dimensions)
the rank $r$ can be choosen smaller than the number of PDEs in an
(off-grid) lattice Boltzmann method. Then, from that perspective,
the amount of memory required and the number of memory operations
we have to perform is reduced. On the other hand, the number of arithmetic
operations is increased. This is precisely the kind of numerical algorithm
that is expected to perform very well on the next generation of supercomputers
(i.e.~exascale systems). Also such algorithms are desperately needed
to fully exploit accelerators, such as graphic processing units and
the Intel Xeon Phi. For more information we refer the reader to the
ASCAC report on exascale computing \cite{ascac}.

\section{The numerical algorithm in the inviscid limit\label{sec:algorithm-in-the-limit}}

An important consideration for the present algorithm is the limit
$\epsilon\to0$. As has been outlined in section \ref{sec:Boltztmann-to-fluid},
the continuous problem (i.e.~the Boltzmann equation) converges to
the Euler equations in this case. To put this statement in the present
framework, in the limit $\epsilon\to0$ the solution of the Boltzmann
equation yields a Maxwell\textendash Boltzmann distribution in phase
space.

In general, however, there is no guarantee that a numerical approximation
conserves this behavior. However, in the present section we show that
each part of the projector-splitting satisfies a very similar constraint.

First, we consider the evolution equations for $K_{j}$ (i.e.~equation
(\ref{eq:evolution-K})). If we take $\epsilon\to0$ we obtain the
constraint
\[
K_{j}-c_{j}^{3}(K)(x)\rho(K)(x)=0
\]
which can be written as
\[
K_{j}=\frac{\rho(K)}{(2\pi)^{d/2}}\int V_{j}(v)\exp\left(-\tfrac{1}{2}(v-u(K))^{2}\right)\,\mathrm{d}v.
\]
This is just the projection of the Maxwell\textendash Boltzmann distribution
onto the space spanned by the $V_{j}$. Thus, as long as $\exp\left(-\tfrac{1}{2}(v-u(K))^{2}\right)$
can be represented accurately in the low-rank manifold (which as we
have discussed in section \ref{sec:Boltztmann-to-fluid} is indeed
the case for weakly compressible flows) this sub-flow of the splitting
algorithm respects the constraints imposed by the continuous problem.

Now, let us consider the evolution equations for $L_{i}$ (i.e.~equation
(\ref{eq:evolution-L})). For $\epsilon\to0$ we obtain
\[
L_{i}=\int X_{i}\rho(L)\exp\left(-\tfrac{1}{2}(v-u(L))^{2}\right)\,\mathrm{d}x.
\]
This takes the Maxwell\textendash Boltzmann distribution and projects
it onto the space spanned by the $X_{i}$. Thus, we again conclude
that if the corresponding low-rank manifold can accurately represent
the Boltzmann\textendash Maxwell distribution our numerical algorithm
will naturally enforce the corresponding constraint. 

Having considered both the evolution equations for $K_{j}$ and $L_{i}$,
it should come as no surprise that we obtain a very similar result
for the evolution equations for $S_{ij}$ (i.e.~equation (\ref{eq:evolution-S})).
In this setting we obtain
\[
S_{ij}=\int X_{i}V_{j}\rho(S)\exp\left(-\tfrac{1}{2}(v-u(S))^{2}\right)\,\mathrm{d}(x,y)
\]
 which once again is just the projection onto the space spanned by
the $X_{i}V_{j}$. Thus, if we can assume that our low-rank approximation
is able to exactly represent the Maxwell\textendash Boltzmann distribution,
then the dynamical low-rank splitting would yield exactly the correct
form of the distribution function $f$.

\section{Numerical results\label{sec:numerical-results}}

In this section we will perform numerical simulations with the proposed
algorithm. As a comparison we consider a classic fluid solver that
uses the second order MacCormack method. 

\subsection{Propagation of sound waves}

As the first test case we consider a simple plane wave propagating
in the $y$-direction. It can be easily shown that, if we can neglect
the nonlinear term in the Navier\textendash Stokes equations (i.e.~for
small velocities), the damped wave equation
\[
\partial_{tt}\rho+\mu\Delta(\partial_{t}\rho)=\Delta\rho
\]
is obtained. Due to the ideal gas law $p=\rho$ this can also be written
as a pressure wave. For small damping (i.e.~small viscosity $\mu$)
we obtain the plane wave solutions
\begin{equation}
\rho(t,x,y)=1+\delta\sin(k_{x}x-\omega t)+\delta\sin(k_{y}y-\omega t)\label{eq:plane-wave-rho}
\end{equation}
with $\omega$ the frequency, $(k_{x},k_{y})$ the wave vector, and
$\delta$ the amplitude of the wave. Since the speed of sound is equal
to unity, frequency and wave vector are coupled by the dispersion
relation $\omega^{2}=k_{x}^{2}+k_{y}^{2}$.

For the numerical example we consider the initial value
\begin{align*}
\rho(0,x,y) & =1+\delta\sin(2\pi y)\\
u_{1}(0,x,y) & =0\\
u_{2}(0,x,y) & =\delta\sin(2\pi y)
\end{align*}
on the domain $\Omega=[0,1]\times[0,1]$. As described above, for
small viscosity and $\delta\ll1$ this results in a plane wave solution
traveling in the $y$-direction with unit speed. An interesting point
here is that the step size of any explicit numerical method would
be dictated by the CFL condition imposed by the speed of sound. That
is, it would have to satisfy $\tau\leq h$, where $\tau$ is the time
step size and $h$ is the grid spacing. On the other hand, the dynamic
low-rank splitting proposed here should be able to take time steps
dictated by accuracy (i.e.~time steps that are significantly larger).

\begin{figure}
\begin{centering}
\includegraphics[width=12cm]{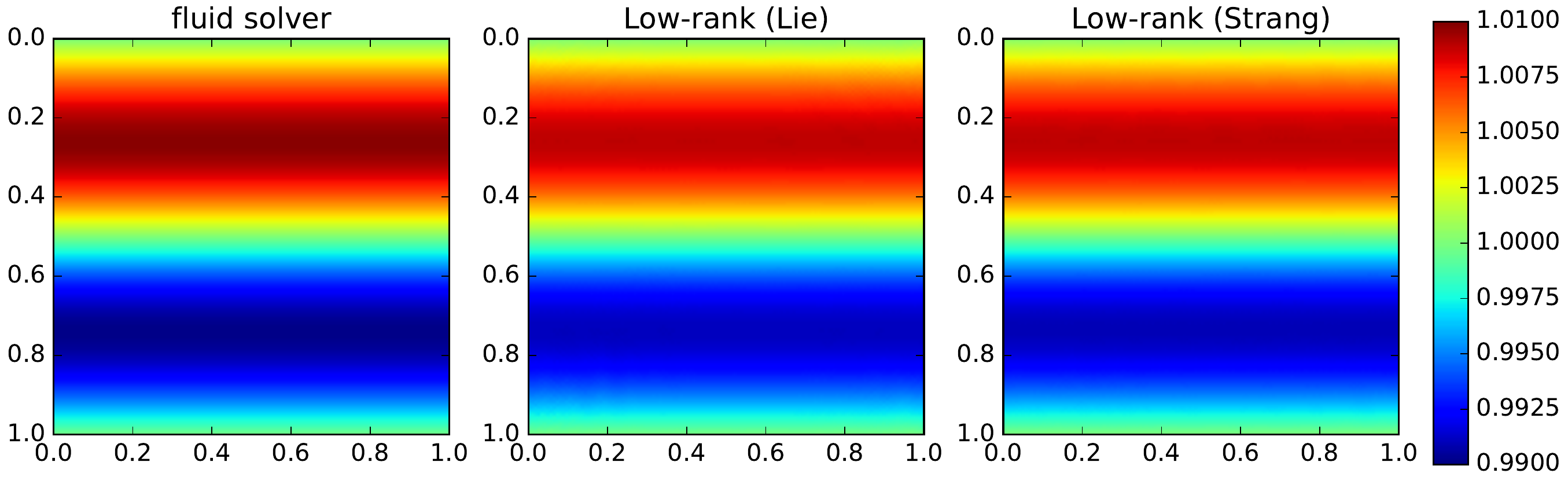}
\par\end{centering}
\caption{The density $\rho$ at time $t=1$ is shown for the classic fluid
solver, Lie splitting, and Strang splitting (with $\epsilon=10^{-3}$).
For all simulations centered differences with $128$ grid points per
direction are used and the rank is set to $10$. The time step size
for the classic fluid solver is set to $\tau=7\cdot10^{-3}$ (a CFL
number of $0.9$). For Lie splitting $\tau=0.1$ and for Strang splitting
$\tau=0.2$ has been used.\label{fig:sound}}
\end{figure}

The results presented in Figure \ref{fig:sound} are meant to check
this reasoning and to verify the code in this simple setting. We observe
that with the dynamical low-rank Strang splitting we can take almost
$30$ times larger time steps compared to the MacCormack method. The
low rank approximation does not conserve mass exactly. However, in
this setting the conservation of mass is still acceptable (on the
order of $10^{-6}$), especially considering that we take quite large
time steps.

\subsection{Shear flow}

Here we consider a shear flow that is given by
\begin{align}
\rho(0,x,y) & =1\nonumber \\
u_{1}(0,x,y) & =v_{0}\begin{cases}
\tanh\left(\frac{y-\tfrac{1}{4}}{\Delta}\right) & y\leq\tfrac{1}{2}\\
\tanh\left(\frac{\tfrac{3}{4}-y}{\Delta}\right) & y>\tfrac{1}{2}
\end{cases}\label{eq:shear-flow-u0}\\
u_{2}(0,x,y) & =\delta\sin(2\pi x).\nonumber 
\end{align}
That is, we have a velocity profile in the $y$-direction that changes
relatively abruptly from $v_{0}=0.1$ to $-v_{0}$ (as we have chosen
$\Delta=1/30$). A small perturbation ($\delta=5\cdot10^{-3}=0.05\cdot v_{0}$)
is then added to the velocity in the $y$-direction. This problem
has been used as a test problem for (mostly incompressible) flow in
a number of publications \cite{bell1989second,liu2000high,einkemmer2014}.

\begin{figure}
\begin{centering}
Fluid solver
\par\end{centering}
\begin{centering}
\includegraphics[width=12cm]{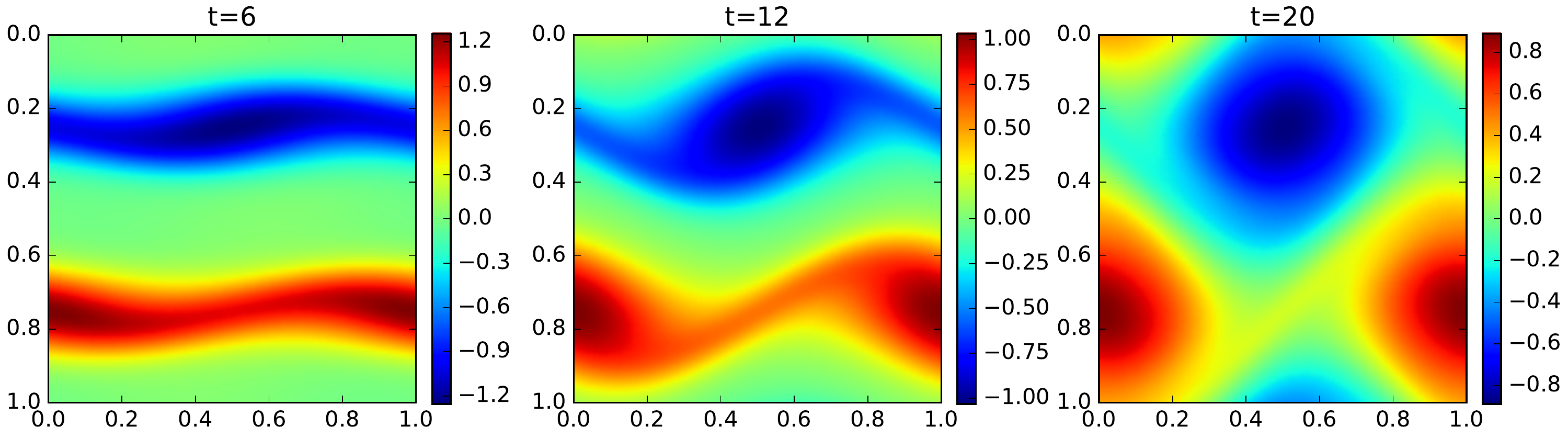}
\par\end{centering}
\begin{centering}
Low-Rank (Strang)
\par\end{centering}
\begin{centering}
\includegraphics[width=12cm]{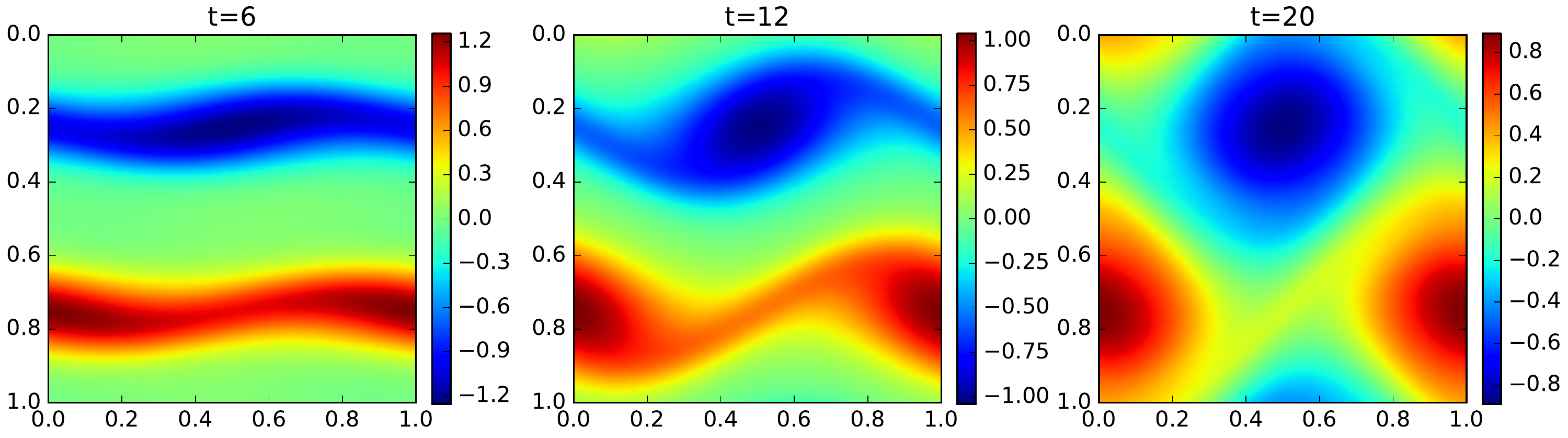}
\par\end{centering}
\centering{}\caption{The time evolution of the vorticity $\omega=\partial_{x}u_{2}-\partial_{y}u_{1}$
for the shear flow given by initial values (\ref{eq:shear-flow-u0})
is shown (with $\text{Re}=300$). The results from a classic fluid
solver are shown on the top and the results from the dynamical low-rank
splitting are shown on the bottom.  For all simulations centered
differences with $128$ grid points per direction are used and the
rank is set to $10$. The time step size for the classic fluid solver
is set to $\tau=7\cdot10^{-3}$ (a CFL number of $0.9$), while for
the low-rank solver the time step size is set to $\tau=0.1$. For
the low-rank implementation Strang splitting is used. \label{fig:shear-Re300}}
\end{figure}

First, we consider a modest Reynolds number ($\text{Re}=300$). The
corresponding results are shown in Figure \ref{fig:shear-Re300}.
As is common for such studies we have plotted the vorticity. We observe
excellent agreement between the proposed low-rank algorithm and the
classic fluid solver. Let us also note that the low-rank algorithm
is not encumbered by a CFL condition. In fact, we can take a time
step that is almost $15$ times as large compared to the fluid solver.

\begin{figure}
\begin{centering}
Fluid solver
\par\end{centering}
\begin{centering}
\includegraphics[width=12cm]{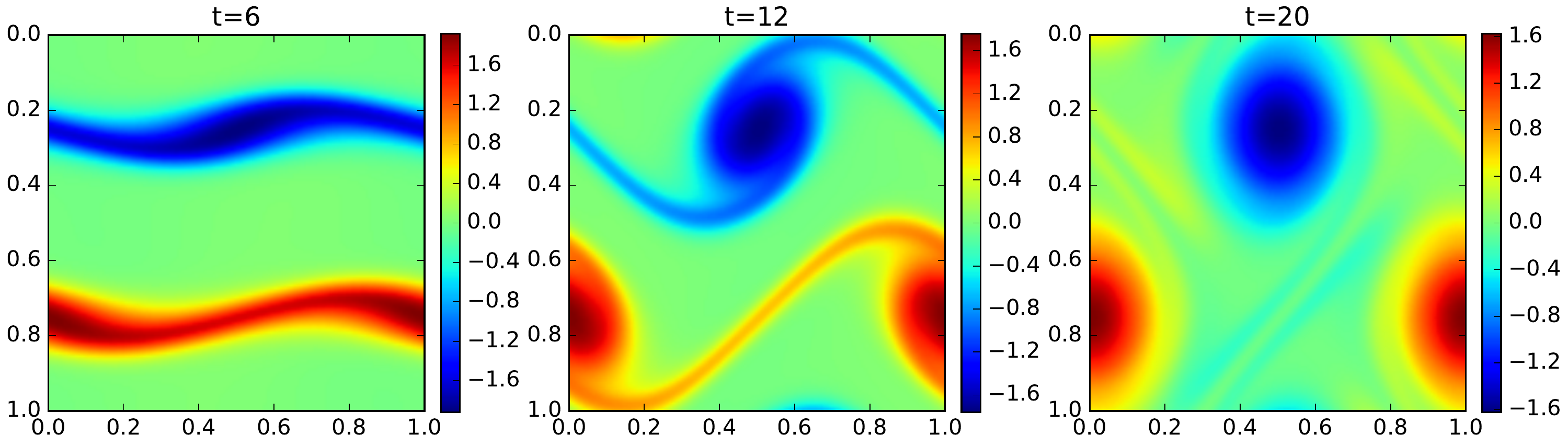}
\par\end{centering}
\begin{centering}
Low-rank (Strang)
\par\end{centering}
\begin{centering}
\includegraphics[width=12cm]{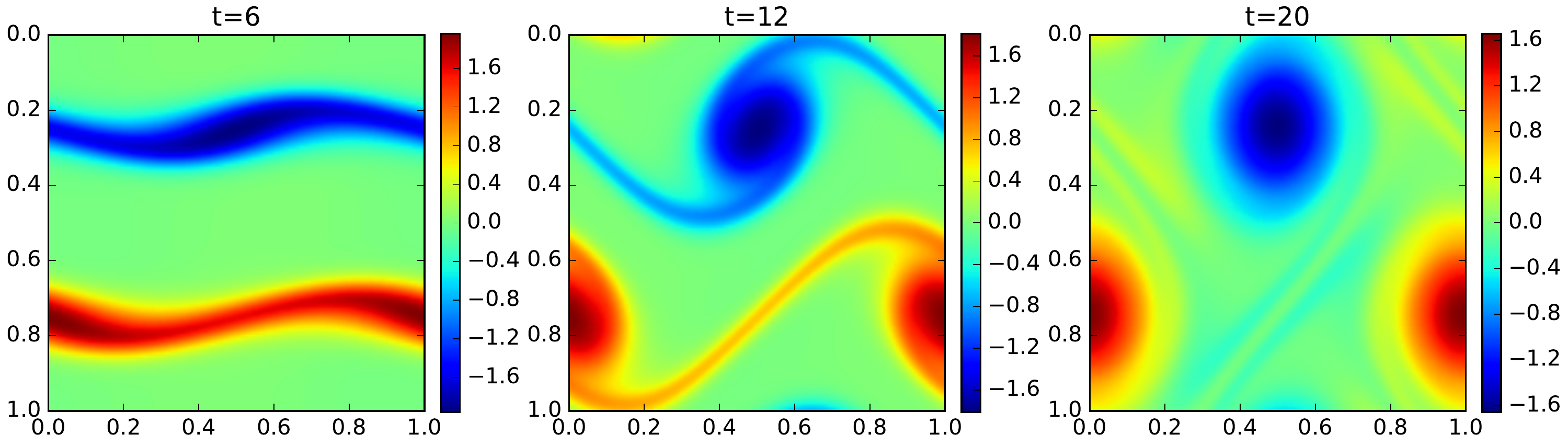}
\par\end{centering}
\centering{}\caption{The time evolution of the vorticity $\omega=\partial_{x}u_{2}-\partial_{y}u_{1}$
for the shear flow given by initial values (\ref{eq:shear-flow-u0})
is shown (with $\text{Re}=1000$). The results from a classic fluid
solver are shown on the top and the results from the dynamical low-rank
splitting are shown on the bottom. For all simulations centered differences
with $128$ grid points per direction are used and the rank is set
to $10$. The time step size for the classic fluid solver is set to
$\tau=7\cdot10^{-3}$ (a CFL number of $0.9$), while for the low-rank
solver it is set to $\tau=0.1$. For the low-rank implementation Strang
splitting is used.\label{fig:shear-Re1000}}
\end{figure}

Second, we increase the Reynolds number to $\text{Re}=1000$. This
is a more challenging problem in the sense that finer structures appear
in the solution. The numerical results are shown in Figure \ref{fig:shear-Re1000}.
We once again observe excellent agreement between our low-rank algorithm
and the classic fluid solver. In fact, all of the conclusions drawn
for the case $\text{Re}=300$ can be applied to the present case as
well.

\begin{figure}
\begin{centering}
Low-rank (Strang)
\par\end{centering}
\begin{centering}
\includegraphics[width=12cm]{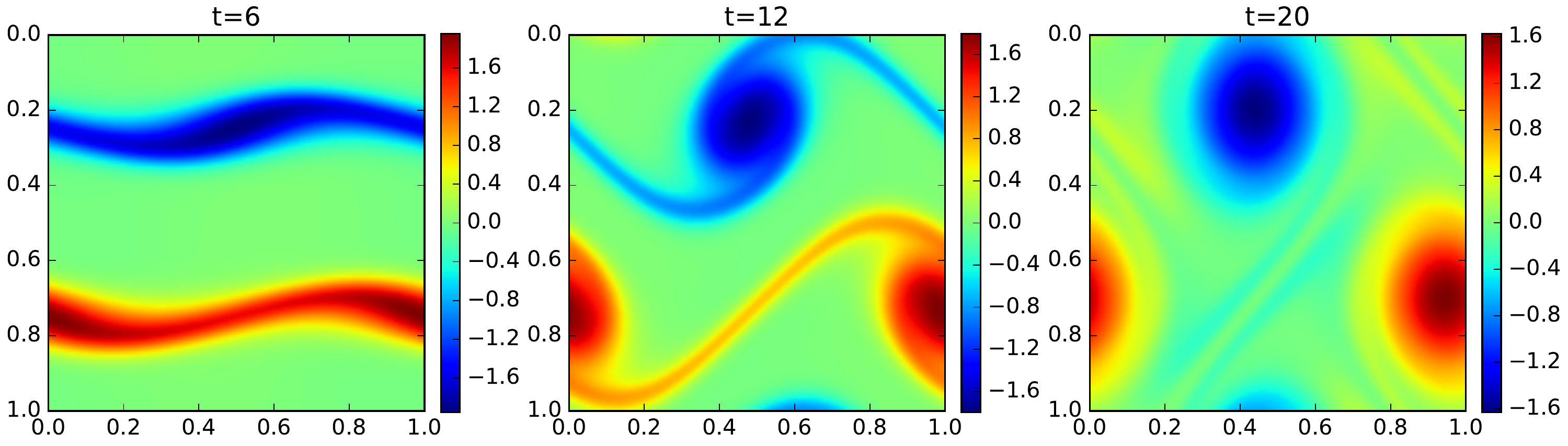}
\par\end{centering}
\centering{}\caption{The time evolution of the vorticity $\omega=\partial_{x}u_{2}-\partial_{y}u_{1}$
for the shear flow given by initial values (\ref{eq:shear-flow-u0})
is shown (with $\text{Re}=1000$). The dynamical low-rank splitting
is used as the integrator. For all simulations centered differences
with $128$ grid points in the space directions and $16$ grid points
in the velocity directions are used. The time step size is set to
$\tau=0.1$ and Strang splitting has been employed.\label{fig:shear-Re1000-nv16}}
\end{figure}

The last point we want to make here is that it is usually not necessary
to use a large number of grid points in the velocity direction. To
demonstrate this, we have repeated our numerical experiment with only
$16$ grid points in the $v$-directions, while still using $128$
grid points in the space directions. In that setting the computational
performance is completely dictated by solving equation (\ref{eq:evolution-K}).
Nevertheless, as Figure \ref{fig:shear-Re1000-nv16} demonstrates,
the numerical results show excellent agreement compared to Figure
\ref{fig:shear-Re1000}, where $128$ grid points where used in the
velocity directions.

\section{Conclusion \& Outlook\label{sec:conclusion}}

We have introduced a numerical algorithm for solving the weakly compressible
Navier\textendash Stokes equations that is based on a dynamical low-rank
splitting algorithm. The behavior of this algorithm has been investigated
and numerical simulations have been conducted that show excellent
agreement with a classic fluid solver.

The algorithm has been considered in the context of weakly compressible
isothermal flow with periodic boundary conditions. However, this restrictions
are not fundamental problems. For example, the extension to temperature
dependent flows is immediate. In fact, only a time and space dependent
$\theta$ has to be introduced in section \ref{sec:Boltztmann-to-fluid}.
This (slightly) changes the collision operator, but the numerical
method remains virtually unaffected. The expansion (\ref{eq:feq-expansion})
is only valid for small velocities (i.e.~weakly compressible flow).
However, this does not mean that we can not efficiently represent
the solution by a low-rank function. In fact, it is not even clear
that equation (\ref{eq:feq-expansion}) is the best low-rank approximation
(i.e.~the approximation with the smallest rank) one can obtain. We
have only considered periodic boundary conditions here. However, similar
to the lattice Boltzmann method, no-slip boundary conditions can be
imposed by a 'bounce-back' scheme. All of this is the subject of future
research. 

Furthermore, the method proposed here offers a path forward for simulations
that need to resolve some kinetic effects. Such problems are common
in various fields of plasma physics. Full scale simulations with the
Boltzmann (collisional Vlasov) equation are often prohibitive from
a computational point of view. However, as has been shown in \cite{einkemmer2018low}
low-rank approximations are still able to resolve a range of kinetic
effects quite well. The method proposed here would thus conceivably
allow us to extend fluid models (say magnetohydrodynamics) to a regime
in which kinetic effects are needed.

\section*{Acknowledgments}

We would like to thank Christian Lubich (University of T\"ubingen)
for the many helpful discussions.

\bibliographystyle{plain}
\bibliography{fluid-lowrank}

\end{document}